\documentstyle[11pt,amssymb,amstex]{article}

\let\goth\frak
\newcommand{\C}{{\Bbb C}}
\newcommand{\ab}{{\hskip 8mm}}
\newcommand{\cc}{{\cal C}}

\newcommand{\oo}{\cal O_}
\newcommand{\J}{{\cal J_q}}
\newcommand{\Hk}{{\bold H_q}}
\newcommand{\HH}{{\bold H}}
\newcommand{\HHH}{{\tilde {\bold H}}}
\newcommand{\hhh}{{\tilde {\cal H}}}
\newcommand{\hk}{{\cal H_q}}

\newcommand{\hdd}{{\vphantom{H}\smash{\buildrel_{\,\bullet\bullet}\over H}}}
\newcommand{\hd}{{\vphantom{H}\smash{\buildrel_{\,\,\bullet}\over H}}}
\newcommand{\xx}{{\bold X}}
\newcommand{\Z}{{\Bbb Z}}
\newcommand{\D}{D}
\newcommand{\g}{{\goth g}}
\newcommand{\h}{{\goth h}}
\newcommand{\dd}{{\Delta}}

\begin{document}
\setlength{\parindent}{0pt}
\setlength{\parskip}{3pt plus 5pt minus 0pt}
\hfill {\small alg-geom/9512017}\break

\centerline {\large{\bf RESIDUE CONSTRUCTION OF HECKE ALGEBRAS}}

\vskip .7cm

\centerline {\sc Victor Ginzburg, Mikhail Kapranov, Eric Vasserot}

\vskip 1cm

\ab\parbox[t]{110mm}{\small Hecke algebras are usually
defined algebraically, via generators and relations.
We give a new algebro-geometric construction of
affine and double-affine Hecke algebras
(the former is known as the Iwahori-Hecke algebra, and
the latter was introduced by Cherednik [Ch1]) in terms of
residues.
More generally, to any generalized Cartan matrix $A$ and a
point $q$ in a 1-dimensional complex
algebraic group $\cc$ we associate
an associative algebra $\Hk$.
If $A$ is of finite type and $\cc=\C^\times$, the algebra
$\Hk$ is the
{\it affine}  Hecke algebra of the corresponding
finite root system. If $A$ is of affine type and $\cc=\C^\times$ then
$\,\Hk$ is, essentially, the Cherednik algebra.
The case $\cc=\C^+$ corresponds to `degenerate' counterparts of the
above objects cosidered by Drinfeld [Dr] and Lusztig [L2].
Finally, taking $\cc$ to be an elliptic curve one gets
some new
elliptic analogues of the affine Hecke algebra.}
\vskip .7cm

\ab Let $W$ be the Weyl group and $T$ the maximal torus
of a Kac-Moody
group associated to the Cartan matrix $A$.
Write $X_*(T)$ for the lattice of one-parameter subgroups
in $T$ and
$\C(T)$
for the field
of rational functions on $T$ with the natural $W$-action.
There is a well-known description of affine Hecke algebra as the
subalgebra of $End_{_\C} \C(T)$ generated by
the so-called  Demazure-Lusztig operators, cf. [Lu].
We observe that all the operators obtained in this way
have poles of a very special kind.
Thus, for any 1-dimensional algebraic group $\cc$,
 we can give an `external' definition of an algebra
$\Hk$ as the subalgebra of operators on the vector space
of rational functions on $\cc\otimes_{_\Z} X_*(T)$
consisting of the elements with a certain
specific singularity type.
All  standard properties
of Hecke algebras (e.g. the normal form of elements of the Cherednik algebra)
follow readily from this definition.  Our construction
completes and somewhat clarifies the results
of Kostant-Kumar [KK1]
on the equivariant cohomology of the affine flag manifold.
\vskip .6cm

{\bf 1. Main construction.}

{\bf 1.1.}$\enspace$
Let $A = \|a_{ij}\|_{i,j=1}^n$
be a symmetrizable generalized Cartan matrix. Throughout the paper we
fix a root data, see [Sp], associated to $A$, that is:

$\bullet\enspace$ A free abelian group $\xx$ such that $rk(\xx) + rk(A)= 2n$;

$\bullet\enspace$ $\Z$-independent elements
$\alpha_1, ..., \alpha_n \in {\xx}$,
called  simple roots;

$\bullet\enspace$ \parbox[t]{115mm}{$\Z$-independent elements
$\alpha^\vee_1, ..., \alpha^\vee_n \in \xx^\vee$, such that
$\langle\alpha^\vee_i, \alpha_j\rangle = a_{ij}$,\\
where $\xx^\vee := Hom(\xx,\Z)$ and
$\langle\,,\,\rangle :$ $\xx^\vee \times \xx\to \Z$ is
the natural pairing.
The $\alpha_i^\vee$
are called simple coroots.}

The vector space $\C\otimes_{_\Z} \xx$ together with the
collections of
simple roots and coroots, as above, is called {\it a realization}
of $A$, see [Ka]. Given a realization of $A$, one associates (see [Ka])
the Kac-Moody Lie algebra $\g(A)$ with Cartan subalgebra
$\C\otimes_{_\Z}\xx^\vee$.
Similarly, given a root data, one associates the Kac-Moody group
$G(A)$ with maximal torus $\C^*\otimes_{_\Z} \xx^\vee$ (the group $G(A)$
will not be used in this paper).

\ab Write $R\subset \xx$ for the set of roots
and $W$ for the Weyl group of $G(A)$.
The group $W$ is generated by the simple reflections $s_i$ , $i=1, ..., n$,
corresponding to the $\alpha_i$. Let
$R^{re}\subset R$ be the set of real roots, i.e., roots of the form
$w\alpha_i$, $w\in W$, $i=1,..., n$. For any $\alpha\in R^{re}$
the corresponding reflection, $s_\alpha$, belongs to
$W$, in particular, $s_i = s_{\alpha_i}$.
One has a decomposition $R=  R_+\sqcup R_-$ into positive
and negative roots, and we put $R^{re}_\pm :=
R^{re}\cap R_\pm$.

\ab From now on, we fix $\cc$, a 1-dimensional complex algebraic group written
multiplicatively. We form
the abelian algebraic group
$T = \cc\otimes_{\Z} \xx^\vee$, and write
$\C(T)$ for the field of {\it rational} functions
on $T$.
Associated to any
$\lambda\in \xx$ is a  group homomorphism $\cc\otimes_{\Z} \xx^\vee \to
\cc$ denoted $t \mapsto t^\lambda$.
Explicitly, if $t=c\otimes x$, where
$c\in \cc, x\in\xx^\vee$, we have $t^\lambda= c^{\langle\lambda,x\rangle}$.
Given a point $q \in \cc$ we put $T_{\lambda,q} :=
\{t^\lambda = q^{}\}$, a divisor in $T$.
This is a reduced but not necessarily connected subvariety of $T$.
If $q=1$ we simply write $T_{\lambda}=T_{\lambda,1}$
for the kernel divisor.

{\bf 1.2.}$\enspace $There is a  $W$-action on $T = \cc\otimes_{\Z} \xx^\vee$
induced by the natural  $W$-action on the second factor.
This gives a left $W$-action
$w : f \mapsto {^w\!\!f}$ on the field ${\C}(T)$
by the formula ${^w\!\!f}(t) = f(w^{-1}\cdot t)$.
Introduce a {\it twisted group algebra},
${\C}(T)[W]$, as the complex vector space
${\C}(T)\otimes_\C \C[W]$ with multiplication:
\[(f\otimes w) \cdot (g\otimes y) = (f\cdot {^w\!\!g})\otimes (w\cdot
y)\quad f,g \in {\C}(T),\; w,y \in W\]
To simplify notation, we will write $f [w]$ instead of
$f\otimes w$ in the future.

\ab Throughout the paper we fix an element $q \in \cc$, $q\neq \pm 1$.
Also, the curve $\cc$ is assumed to be affine, with the exception of
$\S 4$ devoted to the elliptic case.
\hfill\break

{\bf Definition 1.3.} Let $\HHH$ be the
be the ${\C}$-linear subspace in ${\C}(T)[W]$ consisting of
the elements $\sum_{w\in W} f_w [w]$ such that:

\vskip .2cm

(1.3.1) \parbox[t]{115mm}{
Each function $f_w$ has no other singularities but first order poles
at the divisors $T^\alpha$, for a finite number of roots
$\alpha\in R^{re}_+$.}

(1.3.2) \parbox[t]{115mm}{
For every $w\in W$ and $\alpha\in R^{re}_+$, we have
${\rm Res}_{_{T_\alpha}} (f_w) +{\rm Res}_{_{T_\alpha}} (f_{s_\alpha w})
=0.$}

Let $\Hk$ be the subspace in $\HHH$ consisting of elements as above satisfying
the following additional condition:

(1.3.3) \parbox[t]{110mm}{The function
$f_w$ vanishes on $T_{\alpha,q^{-2}}$
whenever $\alpha\in R^{re}_+$ and $w^{-1}(\alpha)\in  R_-$.}

\vskip .3cm
{\bf Theorem 1.4.}  {\it  Both $\HHH$ and $\Hk$
are subalgebras
in ${\C}(T)[W]$.}
\vskip .3cm

{\bf 1.5.} {\sc Proof :}  We first show that
$\HHH$ is a subalgebra.

 Fix a root $\alpha\in R^{re}_+$.
Let
$P = \sum P_w [w]$ and $Q =\sum Q_w [w]$ be two elements
 of $\HHH$. Then
$$P\cdot Q = \sum_{u\in W} F_u [u],\quad\mbox{where}\quad
F_u = \sum_{wy=u}
P_w\cdot ^w\!\!Q_y
\eqno (1.5.1)$$
To check (1.3.1) we must
show that each coefficient $F_u$ has at most first
order pole at  generic points of the hypersurface $T_\alpha$.
To that end, observe that, for fixed $u\in W$,
 the set of pairs
$(w,y)$ such that $wy =u$ has
a fixpoint free involution $(w, y) \mapsto
(s_\alpha w, w^{-1} s_\alpha w y)$.
Split the sum in (1.5.1) into partial sums, each
taken over an orbit of this involution. A
partial sum has the form
$$P_w\cdot {^w\!Q}_y +P_{s_\alpha w}
\cdot {^{s_\alpha w}\!Q}_{w^{-1} s_\alpha wy}
\eqno (1.5.2)$$
Adding and subtracting the term $P_{s_\alpha w}\cdot{^w\!Q}_y$ and using
that
 $w^{-1} s_\alpha w = s_{w^{-1}(\alpha)}$, one rewrites (1.5.2)  as follows
$$\left(P_w+P_{s_\alpha w}\right)\cdot {^w\!Q}_y \,-\,
P_{s_\alpha w}\cdot\left({^w\!Q}_y -{^{s_\alpha
w}\!Q}_{s_{w^{-1}(\alpha)}y}\right)
\eqno (1.5.3)$$
Observe now that the function $P_w+P_{s_\alpha w}$ is regular at
$T_\alpha$ by
(1.3.2), while the function $^w\!Q_y$ has at most simple pole at $T_\alpha$.
Hence,
the first summand in (1.5.3) has at most simple pole at $T_\alpha$.
Similarly, condition (1.3.2) applied to $Q_y$ and the
root $w^{-1}(\alpha)$  says that the function
$Q_y +Q_{s_{w^{-1}(\alpha)}y}$ is regular at the divisor $T_{w^{-1}(\alpha)} =
w^{-1}(T_\alpha)$.
Hence, the function ${^w\!Q}_y +{^{w}\!Q}_{s_{w^{-1}(\alpha)}y}$
which takes an element $t\in T$ into
$Q_y(w^{-1}t) + Q_{s_{w^{-1}(\alpha)}y}(w^{-1}t)$,
is regular for  $t\in T_\alpha$.
Observe finally that, for any function $f$, we have
${\rm Res}_{_{T_\alpha}} (f) = -{\rm
Res}_{_{T_\alpha}}({^{s_\alpha}\!\!f})$.
Therefore the expression
${^w\!Q}_y -{^{s_\alpha w}\!Q}_{s_{w^{-1}(\alpha)}y}$
in the second term of (1.5.3)
is regular at $T_\alpha$. Thus, (1.5.3) has at most simple pole at
$T_\alpha$,
and condition (1.3.1) holds for $P\cdot Q$.

\ab We now check that condition (1.3.2) holds for $P\cdot Q$ and any root
$\alpha$.
Given $u\in W$, we must show that the sum (cf. (1.5.1)):
$$F_u +F_{s_\alpha u}= \sum_{wy=u}
(P_w\cdot {^w\!Q}_y +P_{s_\alpha w} \cdot {^{s_\alpha w}\!}Q_y)
\eqno (1.5.4)$$
is regular at $T_\alpha$.
To that end,
we split the sum in (1.5.4) into partial summs
 corresponding to $w$ running over an individual
orbit of the involution $w\mapsto s_\alpha w$ and $y$ running over an
individual orbit of the involution $y\mapsto w^{-1}s_\alpha w y$
respectively.
Such a partial sum has the form
$$P_w\cdot {^w\!Q}_y +P_{s_\alpha w}
\cdot {^{s_\alpha w}\!Q}_{s_{w^{-1}(\alpha)}y} +
P_{s_\alpha w}\cdot {^{s_\alpha w}\!Q}_y+
P_w\cdot {^{s_\alpha w}\!Q}_{s_{w^{-1}(\alpha)}y}$$
Regrouping terms, we obtain
$$P_w \cdot ({^w\!Q}_y +
{^w\!Q}_{s_{w^{-1}(\alpha)}y})
+ P_{s_\alpha w}\cdot {^{s_\alpha}\!(}Q_{s_{w^{-1}(\alpha)}y}
+{^{w}\!Q}_y)\eqno (1.5.5)
$$
Applying condition (1.3.2) we see as above that in (1.5.5) the two expressions
in the brackets are both regular at $T_\alpha$. Furthermore, observe that
the restrictions of  ${^w\!Q}_y +
{^w\!Q}_{s_{w^{-1}(\alpha)}y}$ and
${^{s_\alpha}\!(}Q_{s_{w^{-1}(\alpha)}y}
+{^{w}\!Q}_y)$ to $T_\alpha$ are the same. Writing $r$ for this
restriction
we find that
the residue of (1.5.5) at $T_\alpha$ equals
\[r\cdot {\rm Res}_{T_\alpha}(P_w +P_{s_\alpha w})\]
The second factor vanishes by (1.3.2).
Thus, we have proved that $\HHH$ is a subalgebra.

{\bf (1.5.6) Notation :} given $w \in W$,  write ${\D}(w)$ for
the set of $\alpha\in R^{re}_+$
such that $w^{-1}(\alpha)\in  R_-$.

\ab We now  prove (1.3.3) holds for $P\cdot Q$
provided it holds for both  $P$ and $Q$.
Let $u\in W$. Then $\alpha\in {\D}(u)$ implies
 $\alpha\in R^{re}_+$ and $u^{-1}(\alpha)\in  R_-$.
 We claim that every summand, $F_u$, in (1.5.1) vanishes at generic points of
$T_{\alpha,q^{-2}}$. Indeed, let $wy=u$. Then
$y^{-1}w^{-1}(\alpha)\in  R_-$, and there are two alternatives:
either $w^{-1}(\alpha)\in  R_+$, or $w^{-1}(\alpha)\in  R_-$.
If $w^{-1}(\alpha)\in  R_+$,
then $w^{-1}(\alpha) \in {\D}(y)$, hence ${Q}_y$
vanishes on $T_{w^{-1}(\alpha),q^{-2}}$. This means that the factor
${^w\!Q}_y$ in (1.5.1) vanishes along $T_{\alpha,q^{-2}}$.
Similarly, if $w^{-1}(\alpha)\in  R_-$, then the factor $P_w$ in (1.5.1)
 vanishes on $T_{\alpha,q^{-2}}$.

\ab Thus in both cases the summand in $F_u$ corresponding to any given
$(w,y)$ vanishes on $T_{\alpha,q^{-2}}$, and  the theorem follows.
$\quad\square$
\vskip .3cm

\noindent {\bf (1.6) Remarks.} (a)
We can, if we want, introduce as many
``quantization parameters" $q_\nu$ as there are $W$-orbits on $R^{re}$,
 by modifying (1.3.3) in an obvious way. Proof of
theorem 1.4 still goes through.
\vskip .1cm

(b) One can modify definition 1.3 letting
the point $q$ vary. This way one gets an algebra $\HH_\cc$
over $\C[\cc]$, the ring of regular functions on $\cc$
(if $\cc=\C^\times$ we have $\C[\cc]={\C}[q, q^{-1}]$, in accordance with
standard definition).
\vskip .1cm

(c) If the Cartan matrix $A$
is of affine type the Lie algebra
$\g(A)$ differs from its derived algebra ${\g}_{der}(A) = [{\g(A)},
{\g}(A)]$. It is then possible to take $\xx^\vee$ as a lattice
in $\h_{der}$, the Cartan subalgebra of ${\g}_{der}(A)$ rather than
$\g(A)$. One then regards the roots as elements of the $\h_{der}^*$
(see $\S 6$) and writes $\Phi_{der}$ for the resulting root data.

\vskip .3cm
\ab The following results will be proved in sections 5 and 6 respectively.

\vskip .3cm
{\bf Theorem  1.7.}  {\it If $A$ is of finite type, then $\Hk$
is isomorphic to the affine Iwahori-Hecke algebra associated to
$A^t$, the transpose of $A$.}

\vskip .3cm
{\bf Theorem 1.8.} {\it If $A$ is of  affine type, then
the algebra $\Hk$ associated to
the root data $\Phi_{der}$, is isomorphic to
the  double affine Hecke algebra (cf. \S6).}

\vskip .3cm

\pagebreak[3]
{\bf 2.  Geometric interpretation.}

\ab We have found  somewhat `strange' conditions (1.3.1-3)
as a result of computation of the equivariant K-theory of
the Steinberg variety. The latter is related to Hecke algebras
due to the work of Kazhdan-Lusztig [KL] (in the finite case);
see Garland-Grojnowski announcement [GG] in the general case.
 An alternative
geometric interpretation of
the algebra $\Hk$ will be given in
theorem 2.2 below.

\ab Assume first that
$A$ is a Cartan matrix of finite
type. Then we
have $rk(A)=n=\dim T$, and $W$ is a finite reflection group.
The algebra $\C[T]$ of regular functions on $T$ is,
 by the Pittie-Steinberg theorem,
a free module over $\C[T]^W$, the
subalgebra of $W$-invariants.
Observe next that
there is a natural  $\C[T]^W$-linear action of the algebra
$\C(T)[W]$ on $\C(T)$. The action of $P = \sum P_w [w]$
is given
by the formula
$$f \mapsto {\hat P}(f) :=\sum_{w\in W} P_w\cdot{^w\!\!f} \eqno(2.1)$$

\ab Let $J_q\subset\C[T]$ denote the principal ideal
generated by the element
$$\dd=\prod_{\alpha\in  R_+} (q^{-1}t^{\alpha/2} - qt^{-\alpha/2}) \in \C[T]$$
\vskip .1cm

{\bf Theorem 2.2.} {\it If $A$ is a Cartan matrix of finite
type, then
the assignment $P\mapsto {\hat P}$, see (2.1), yields
a natural algebra isomorphism}:
$$
\Hk \simeq \{ u\in End_{_{\C[T]^W}} \C[T]
 \enspace | \enspace u(J_q) \subset J_q \}.$$

\ab We will prove Theorem 2.2 in the next section. It will be deduced from
 the following

\vskip .1cm

{\bf Proposition 2.3.} {\it If $A$ is of finite type, then
the assignment $P\mapsto {\hat P}$
yields
an algebra isomorphism}
$$\HHH = End_{_{\C[T]^W}} \C[T] \;.$$
The isomorphism shows that $\HHH$
is a matrix algebra over $\C[T]^W$. Observe further that
theorem 2.2 reads : $\Hk\simeq End_{_{\C[T]^W}} \C[T]
\,\cap\,End_{_{\C[T]^W}} J_q$. Using that, $J_q= \Delta\cdot\C[T]$,
we obtain the following
\vskip .1cm

{\bf Corollary 2.4.} $\;\Hk= \HHH\,\cap\,\Delta\cdot\HHH\cdot\Delta^{-1}\quad
\square$.
\vskip .3cm

\ab In the rest of this section we assume
$A$ to be an arbitrary Cartan matrix and
study the structure
of the algebra $\Hk$ in more detail.
\vskip .2cm
\ab
Observe first that the algebra $\C[T]$
may be identified with a (non-central) subalgebra
of $\Hk$ via the assignment $f\mapsto f[1]$. Further,
for each $i=1,2,\ldots,n$,  introduce the
following element:
$$\sigma_i = \biggl( {qt^{\alpha_i}-q^{-1}\over t^{\alpha_i}-1}\biggl) [s_i]
\,\,
- \,\, {q-q^{-1}\over t^{\alpha_i}-1} \, [1] \in {\C}(T)[W].
\eqno (2.4)$$
These elements belong to $\Hk$ and are called the Demazure-Lusztig
operators, cf. [L1]. Moreover, it is well known (see, e.g., [KL])
that  the elements $\sigma_i\,,\,i=1,2,\ldots,n$ satisfy the
braid relations: if $w\in W$ is any element, then for any reduced
decomposition $w=s_{i_1}\cdot\ldots\cdot s_{i_k}$ the product
$\sigma_w=\sigma_{i_1}\cdot\ldots\cdot \sigma_{i_k}$ depends only on
$w$ but not on the choice of a reduced  decomposition.

\vskip .3cm
{\bf Theorem 2.5.} {\it Let $A$ be an arbitrary Cartan matrix.
 Then $\Hk$ is a free left
$\C[T]$-module with the base} $\{\sigma_w\,,\,w\in W\}$.
\vskip .3cm

\vskip 2mm

\ab Let $H^\sigma$ denote
the
left
$\C[T]$-submodule in $\C(T)[W]$ generated by $\{\sigma_w\,,$
$w\in
W\}$.
We will prove theorem 2.5 by  showing
that $\Hk=H^\sigma$.
We proceed by induction on
`$\leq$', the Bruhat (partial) order on the Weyl group $W$.
For each $w\in W$, introduce the subspace
\[\C(T)[W]_{_{\leq w}} = \{ \sum_{y\in W} f_y [y] \,
 \in \C(T)[W]\, |\,\; f_y = 0
\quad {\rm unless}\quad y\leq w \}\]
and similarly define $\C(T)[W]_{_{< w}}$. This gives an increasing
filtration on $\C(T)[W]$
labeled by the partially ordered set $W$. We write
${\bold H}_{_{\leq w}}$ and $H^\sigma_{_{\leq w}}$
for the corresponding induced filtrations on the subspaces
$\Hk$ and $H^\sigma$ respectively.

\ab Further, given $w\in W$, write ${\C}[T]\theta_w$ for the rank one
${\C}[T]$-submodule in $\C(T)$ with generator
$\theta_w=\prod_{\alpha\in {\D}(w)}
\left( {qt^{\alpha}-q^{-1}\over t^{\alpha}-1}\right)$ (notation 1.5.6).
\vskip 1mm

{\bf Lemma 2.8.} (i) {\it
Let
$f=\sum_{y\leq w} f_y [y]\,\in{\bold H}_{_{\leq w}}$.
Then $f_w \in {\C}[T]\theta_w$.}

\ab (ii) $\sigma_w - \theta_w [w] \,\in\,\C(T)[W]_{_{< w}}$.
\vskip 2mm

\noindent {\sc Proof:} (i) Observe that
$s_\alpha\cdot w > w$ whenever $\alpha\in {\D}(w)$, cf. [Bou].
 Hence, for $f\in {\bold H}_{_{\leq w}}$ and
$\alpha \not\in R_+\,\setminus\,D(w)$, we have  $f_{s_\alpha w}=0$, and
conditions
(1.3.1-3) say that
\begin{enumerate}
\item The only possible singularities of $f_w$ are first order poles
at the divisors $T_{\alpha}$ for $\alpha\in {\D}(w)$;
\item $f_w$ vanishes at
$T_{\alpha, q^{-2}}$ for $\alpha\in {\D}(w)$.
\end{enumerate}
The space of rational functions subject to these conditions
is precisely the free ${\C}[T]$-submodule generated by $\theta_w$.
This proves part (i).

\ab To prove claim (ii), write a reduced expression
$w=s_{i_1}\cdot\ldots\cdot s_{i_k}$. Observe that
the
collection of roots
$$\{\alpha_{i_k}\,,\,s_{i_{k-1}}(\alpha_{i_k})\,,\,
s_{i_{k-2}}s_{i_{k-1}}(\alpha_{i_k})\,,\ldots,\,
s_{i_1}\ldots s_{i_{k-1}}(\alpha_{i_k})\}\,=\, D(w)\eqno(2.9)$$
gives a particular enumeration
of all the elements of the set $D(w)$, cf. [Bou].
Write $\theta_\alpha :=  {qt^{\alpha}-q^{-1}\over t^{\alpha}-1}$,
and recall that $\sigma_{i}=\theta_{\alpha_{i}}[s_i] +
{q-q^{-1}\over t^{\alpha}-1}[1]$.
Thus, by definition of $\sigma_w$ we find
\begin{align*}
\sigma_w
& =\sigma_{i_1}\cdot\ldots\cdot \sigma_{i_k}\\
& =\sigma_{i_1}\cdot\ldots\cdot \sigma_{i_{k-1}}\cdot
(\theta_{\alpha_{i_k}}[s_{i_k}]) + \quad\mbox{{\small {\it lower terms}}}
\\
&=\sigma_{i_1}\cdot\ldots\cdot \sigma_{i_{k-2}}\cdot
({^{s_{i_{k-1}}}\!\theta}_{\alpha_{i_k}}[s_{k-1}s_{i_k}])
 + \quad\mbox{{\small {\it lower terms}}} =\\
& \ldots\ldots\ldots\ldots \\
& =\theta_{\alpha_{i_1}}{^{s_{i_1}}\!\theta}_{\alpha_{i_2}}
{^{s_{i_1}s_{i_2}}\!\theta}_{\alpha_{i_3}}\ldots
{^{s_{i_1}\ldots s_{i_{k-1}}}\!\theta}_{\alpha_{i_k}}
\cdot[s_{i_1}\cdot\ldots\cdot s_{i_k}] + \quad\mbox{{\small {\it lower
terms}}}\\
& =\prod_{\alpha\in {\D}(w)}
\left( {qt^{\alpha}-q^{-1}\over t^{\alpha}-1}\right)[w] +
\quad\mbox{{\small {\it lower terms}}}=
\theta_w [w] + \quad\mbox{{\small {\it lower terms,}}}
\end{align*}
where `{\it lower terms}' stands for a linear combination
$\sum_{y<w} P_y [y]$. The lemma follows.$\quad\square$

\ab {\sc Proof of theorem 2.5:} Theorem 1.4 implies that,
for any $y\in W$, we have $\sigma_y \in \Hk$. Clearly,
$H^\sigma_{_{\leq w}}$ is a free left $\C[T]$-module with
base $\{\sigma_y\,,\, y\leq w\}$. Whence, $H^\sigma_{_{\leq w}}
\subseteq \HH_{_{\leq w}}$. We prove that this inclusion is actually
an equality by induction on $w$ (relative to the Bruhat order).

\ab Let $f=\sum_{y\leq w} f_y [y] \in \HH_{_{\leq w}}$.
Lemma 2.8(i) shows that $f_w= P\cdot\theta_w$, for some $P\in \C[T]$.
Furthermore, part (ii) of the same lemma yields
$f-P\cdot\sigma_w \in \HH_{_{< w}}$. But the  induction hypothesis
yields $\HH_{_{< w}}=H^\sigma_{_{< w}}$.
Hence, $f= (f- P\cdot\sigma_w) + P\cdot\sigma_w \in \HH_{_{< w}}
+\HH_{_{\leq w}}$,
and the theorem follows.
$\quad\square$

{\bf Remark 2.10.} Theorem 2.2 implies in particular that,
for any simple root $\alpha_i$, one has $\sigma_i(\dd) \in J_q$.
In fact, the following more precise result holds:
\[\sigma_i(\dd)=-q^{-1}\cdot\dd.\]

{\sc Proof:} For a simple root $\alpha_i$, set
$\;\Delta_i=q^{-1}t^{-\alpha_i/2} - qt^{\alpha_i/2}$.
Direct calculation gives
$\Delta_i=(\sigma_i-q)t^{\alpha_i/2}$.
Hence, the Hecke relation $(\sigma_i+q^{-1})(\sigma_i-q)=0$ yields
$$(\sigma_i+q^{-1})\Delta_i=
(\sigma_i+q^{-1})(\sigma_i-q)t^{\alpha_i/2}=0\eqno(2.11)$$
This proves the claim in the ${\goth s}{\goth l}_2$-case.

\ab In general,
fix a simple root $\alpha_i$, and write
$\dd=\Delta'\cdot\Delta_i$ where
$$
\Delta'=\prod_{\alpha\in  R_+,\alpha\neq\alpha_i}
(q^{-1}t^{\alpha/2} - qt^{-\alpha/2})$$
Observe, that the simple reflection $s_i$ acts by permutation of the
set $R_+\,~\setminus~\,~\{\alpha_i\}$. Therefore, we get
${^{s_i}\!\!\Delta'}=\Delta'$. Thus, from formula (2.1)
we calculate

\begin{align*}
\sigma_i(\dd) &=
\biggl( {qt^{\alpha_i}-q^{-1}\over t^{\alpha_i}-1} [s_i] \,\,
- \,\, {q-q^{-1}\over t^{\alpha_i}-1} \, [1]
\biggl)(\Delta'\cdot\Delta_i)\\
&={qt^{\alpha_i}-q^{-1}\over t^{\alpha_i}-1}
\cdot{^{s_i}\!\Delta'}\cdot{^{s_i}\!\Delta_i}
\,- \, {q-q^{-1}\over t^{\alpha_i}-1}\Delta'\cdot\Delta_i\\
&=\Delta'\cdot {qt^{\alpha_i}-q^{-1}\over t^{\alpha_i}-1}
\cdot{^{s_i}\!\Delta_i}\,
- \,\Delta'\cdot{q-q^{-1}\over t^{\alpha_i}-1}\cdot\Delta_i\\
&=\Delta'\cdot \sigma_i(\Delta_i)\quad\mbox{by (2.11)}\\
&=\Delta'\cdot (-q^{-1})\cdot\Delta_i=-q^{-1}\cdot\dd\qquad\square.
\end{align*}
%
{\bf 3. Polynomial representation.}

\ab In this section we will study in more detail the action of
the algebras $\Hk$ and $\HHH$ on the vector space $\C(T)$
given by formula (2.1) and give proofs of Proposition 2.3
and Theorem 2.2.
The curve
$\cc$ is  supposed to be affine, and the matrix $A$ is supposed
to be of finite type.

{\bf Lemma 3.1.} (i) {\it Let $P\in \HHH$. Then, for any regular
function $f\in \C[T]$, we have ${\hat P}(f)\in \C[T]$.}

\ab (ii) {\it Assume
 $f\in \C[T]$ satisfies the following condition: for each
positive root $\alpha\in R_+$, the function $f$ vanishes
at $T_{\alpha,q^{-2}}$. Then, for any  $P\in \Hk$, the function
${\hat P}(f)$ satisfies similar condition.}

\vskip 2mm

\noindent {\sc Proof:} The argument is very similar to the one
used in the proof of theorem 1.4. To prove (i), fix a
root $\alpha\in R_+$ and write $P = \sum P_w [w]$.
Split the sum ${\hat P}(f)=\sum_w P_w\cdot ^w\!\!f$ into partial sums
of the form
\[P_w\cdot {^w\!f} +P_{s_\alpha w}
\cdot {^{s_\alpha w}\!f}\]
Setting $g={^w\!f}$ and regrouping terms, one rewrites this  as follows
$$(P_w+P_{s_\alpha w})\cdot g \,-\,
P_{s_\alpha w}\cdot(g-{^{s_\alpha}\!g})\eqno(3.2)$$
Observe now that the function $P_w+P_{s_\alpha w}$ is regular at
$T_\alpha$ by
(1.3.2), while the function $g$ is regular everywhere. Hence the
the first summand in (3.2) has no singularity at $T_\alpha$.
Similarly, $P_{s_\alpha w}$
has at most simple pole at $T_\alpha$, while the function
$g-{^{s_\alpha}\!g}$ vanishes at $T_\alpha$. Hence, the second
summand in (3.2) has no singularity at $T_\alpha$.
Thus, we have shown that, for any root $\alpha$, the function
${\hat P}(f)$ is regular at generic points of $T_\alpha$.
Therefore, it is regular everywhere, and part (i) follows.
Part (ii) is proved in a similar way, repeating the
verification of condition (1.3.3) in the proof of theorem 1.4.
$\quad\square$.

\vskip 2mm

\ab Recall next that the ring $\C[T]^W$ is the coordinate ring of the
orbi-space $T/W$, an affine algebraic variety.
Write $I_a\subset\C[T]^W$ for the
maximal ideal corresponding to a point $a\in T/W$. Given
a $\C[T]^W$-module $M$, set $M_a:= M/I_a\cdot M$, the geometric fiber
of $M$ at $a$.

\ab Observe that by
lemma 3.1(i) we have a
well-defined map
$$\Phi :\HHH \to End_{_{\C[T]^W}} \C[T] \,,\eqno(3.3)$$
which is, moreover, a morphism of $\C[T]^W$-modules.
Given a point $a\in T/W$, write
$$\Phi_a :\HHH_a \to End_{_{\C[T]^W}} \left(\C[T]_a\right) \eqno(3.4)$$
for the induced map of geometric fibers at $a$.

\ab Let $\pi: T\to T/W$ be the natural projection. We have the following

{\bf Lemma 3.5.} {\it Let $t\in T$
 be such that one of
the following two conditions holds:
\begin{enumerate}
\item For any $\alpha \in R_+$, we have $t\not\in T_\alpha$;
\item $t\in T_\alpha$ for exactly one root  $\alpha \in R_+$.
\end{enumerate}
Then, for $a=\pi(t)$, the map (3.4) is surjective.}

{\sc Proof:} Recall that $T/W$ is a smooth algebraic variety (e.g.,
if $W$ is the
symmetric
group, $T/W$ is the symmetric power of the curve $\cc$).

\ab Consider case (1) first. Then $t$ is regular and
the map $\pi: T\to T/W$ is unramified over $a$. Hence,
the fiber $\pi^{-1}(a)= \{w(t)\,,\, w\in W\}$ is the  $W$-orbit
consisting of $\#W$ elements. This yields a natural
$W$-equivariant isomorphism of vector spaces $\C[T]_a \simeq
\C[W]$, where the line $\C\cdot [w]\subset \C[W]$ is identified with the
geometric fiber of $\C[T]$ at the point $w(t)$.
On the other hand, theorem 2.5 shows that $\HHH_a$ is a
$(\#W)^2$-dimensional
$\C$-algebra generated by the elements $\{[y]\,,\, y\in W\}$ and
$\{f_a\,,\, f\in \C[T]\}$. Furthermore, it follows from definitions that
the
 map $\Phi_a :\HHH_a \to End_{_\C} \C[W]$ arising from
(3.4) sends $[y]$ to the operator of left translation by $w^{-1}$,
and $f_a$ to the diagonal operator acting on the line $\C\cdot
[w]\subset\C[W]$
via multiplication by $f(w^{-1}(t))$, the falue of $f$ at the
point $w(t)$. It is clear further, that the operators of those
two types generate the whole algebra $End_{_\C} \C[W]$. That proves
the lemma in case (1).

\ab Assume now  $t\in T_\alpha$ as in case (2). Then the isotropy
group of $t$ in $W$ is
$\langle s_\alpha\rangle=\{1, s_\alpha\}$,
the two-element subgroup in $W$. Locally near $t$, the projection
$\pi: T\to T/W$ is a two-branch ramified covering
of the form $(z_1,z_2,\ldots,z_r) \mapsto$
$ (z_1^2,z_2,\ldots,z_r)$.
Hence, the
induced  linear map of the cotangent spaces
$d\pi^*: T_a^*(T/W) \to T_t^*(T)$ has 1-codimensional
image. Furthermore, the root $\alpha \in Lie(T)^*$
$\simeq T_t^*(T)$ does not belong to the image of $d\pi^*$.
Thus, there is a canonical direct sum decomposition
$$T_t^*(T) = d\pi^*\bigl(T_a^*(T/W)\bigr) \oplus \C\cdot\alpha
\eqno(3.6)$$
Write $pr_\alpha : T_t^*(T) \to \C\cdot\alpha$ for the second
projection.

\ab To any regular function $f$ on a neighborhood of $t$
we
assign a vector in the two-dimensional vector space
$\C\oplus \C\alpha$ with base $1$ and $\alpha$.
 The assignment is given by the formula
\[ f \mapsto \tau(f)= f(t)\cdot 1 + pr_\alpha(df_t)\,,\]
where $df_t$ stands for the differential of $f$ at $t$.
Further, make $\C\oplus \C\alpha$ into
$\langle s_\alpha\rangle$-module setting
 $s_\alpha(1)=1$ and $s_\alpha(\alpha)=-\alpha$.
Then, $\tau$ becomes  $\langle s_\alpha\rangle$-equivariant map.
 Moreover, decomposition
(3.6) shows that $\tau(f)$ depends only on
$f_a$, the image of $f$ in $\C[T]_a=\C[T]/I_a\cdot\C[T]$.

\ab It is clear that the algebra $End_{_\C} (\C\oplus \C\alpha)$
is generated by the following operators:
$$A= {{1\mapsto \alpha}\choose {\alpha\mapsto 0}}
\quad,\quad B={{1\mapsto 0}\choose {\alpha\mapsto 2}}
\quad,\quad C=Id $$
On the other hand, consider the $\C(T)[W]$-action on $\C(T)$.
It is straightforward to verify that, for any $f\in\C(T)$ regular at $t$,
one has
$$ \tau\left((t^\alpha-1)\cdot f\right)= A(\tau(f))\;\;,
\;\;\tau\left((\frac{1}{t^\alpha-1}-\frac{1}{t^\alpha-1}[s_\alpha])\cdot
f\right)=
B(\tau(f)) \eqno(3.7)$$

\ab Observe next that we have a $W$-equivariant
isomorphism of sets $\pi^{-1}(a)$
$\simeq W/\langle s_\alpha\rangle$.
It follows the the map $\tau$ defined above can be uniquely
extended to a natural $W$-module isomorphism
$$ \C[T]_a \stackrel{\sim}{\to}
Ind_{\langle s_\alpha\rangle}^W (\C\oplus \C\alpha)=
 \C[W]\otimes_{_{\C[\langle s_\alpha\rangle]}} (\C\oplus \C\alpha)
$$
Thus, we may view (3.4) as a homomorphism
$$\Phi_a : \HHH_a \to End_{_\C}
\left(Ind_{\langle s_\alpha\rangle}^W (\C\oplus \C\alpha)\right)
\eqno(3.8)$$
\ab To prove that $\Phi_a$ is surjective, we observe, as in proof
of case (1), that the image of $\Phi_a$ contains left
translations by $W$ as well as all `diagonal' operators,
given by multiplication by a scalar-valued function on
$W/\langle s_\alpha\rangle$. In addition, $\HHH_a$ contains
the element which is equal to  $t^\alpha-1$, resp.
$\frac{1}{t^\alpha-1}-\frac{1}{t^\alpha-1}[s_\alpha]$, at $t$,
 and equal to 0 at all other points of the fiber $\pi^{-1}(a)$.
Equation (3.7) then shows that the image of $\Phi_a$ contains
all operators that coomute with left translation and act
as the operators $A$, resp. $B$, fiberwise. Altogether, all the
listed operators generate the RHS of (3.8). That completes
the proof of case (2) of the lemma. $\quad\square$

\vskip 2mm
{\sc Proof of proposition 2.3:} By the Pittie-Steinberg theorem,
the
object on the RHS of (3.3) is a free $\C[T]^W$-module of rank
$(\#W)^2$.
Similarly, theorem 2.5 implies that the
object on the LHS of (3.3) is a free $\C[T]^W$-module of the same
rank. Hence, proving the proposition suffices it to show
that the map $\Phi$ in
(3.3) is surjective. Furthermore, the modules being free,
surjectivity of $\Phi$ follows from the surjectivity
of the corresponding maps $\Phi_a$, for all points $a$ in the
complement of a codimension two subset of $T/W$. But
the points $\pi(t)$ such that $t$ satisfies
 neither of the conditions of
lemma 3.5 form a codimension two subset. Thus, proposition
2.3 follows from the lemma. $\quad\square$.

 \vskip 2mm
{\sc Proof of theorem 2.2.:} In view of proposition 2.3,
the theorem is clearly equivalent to the
equality, cf. corollary 2.4:
$\,\Hk= \HHH\,\cap\,\Delta\cdot\HHH\cdot\Delta^{-1}\,$.
Theorem 2.5 and remark 2.10 imply that the LHS of the equality
is contained in the RHS.
To prove the opposite inclusion, let
$P=\sum_w P_w[w] \in \HHH$. Then we calculate

\begin{align*}
\Delta^{-1}\cdot P\cdot\Delta& =\sum_w
\Delta^{-1}\cdot P_w \cdot{^w\!\!\Delta}\\
&=\sum_w\prod_{\alpha\in R_+}
(q^{-1}t^{\alpha/2} - qt^{-\alpha/2})^{-1}\cdot P_w \cdot
\prod_{\alpha\in R_+}
(q^{-1}t^{w(\alpha)/2} - qt^{-w(\alpha)/2})\\
&=\sum_w\prod_{\alpha\in R_+}
(q^{-1}t^{\alpha/2} - qt^{-\alpha/2})^{-1}
\prod_{w^{-1}(\alpha)\in R_+}
(q^{-1}t^{\alpha/2} - qt^{-\alpha/2})
\cdot P_w\\
&=\sum_w\prod_{\alpha\in {\D}(w)}
{q^{-1}t^{w(\alpha)/2} - qt^{-w(\alpha)/2}\over q^{-1}t^{\alpha/2} -
qt^{-\alpha/2}}
\cdot P_w
\end{align*}

We see that the assumption
$\Delta^{-1}\cdot P\cdot\Delta\in \HHH$
implies that, for any $w\in W$, the function
$\prod_{\alpha\in {\D}(w)}
\frac{q^{-1}t^{w(\alpha)/2}
- qt^{-w(\alpha)/2}}{q^{-1}t^{\alpha/2} - qt^{-\alpha/2}}
\cdot P_w$ has no singularities. Whence, $P_w$ is divisible
by $\prod_{\alpha\in {\D}(w)}$
$(q^{-1}t^{\alpha/2} - qt^{-\alpha/2})$,
and condition (1.3.3) follows. $\quad\square$

\vskip 7mm
\pagebreak[3]

{\bf 4. Elliptic case.}

\ab In this section we keep the assumption that the Cartan matrix $A$
 is of finite type
but allow the 1-dimensional algebraic group $\cc$ to be arbitrary,
an elliptic curve in particular. It is convenient,
in this more general setting, to use
the language of sheaves. Write
$\oo{_X}$ for the sheaf of regular functions on an algebraic
variety $X$.

\vskip .1cm
\ab As in section 3, let $\pi: T\to T/W$ be the natural projection.
Observe that the group $W$ acts freely on the Zariski-open subset
$T^{reg} := T\,\setminus\,(\cup_{\alpha\in R} T_\alpha)$ so that
$\pi$ is a finite flat morphism
unramified over $\pi(T^{reg})$. It follows from the flatness
that for any locally free sheaf ${\cal O}_T$-sheaf $\cal L$
its direct image, $\pi_* {\cal L}$, is a locally free sheaf on $T/W$.

\ab Let $j: T^{reg}\hookrightarrow T$ be the embedding.
It will be useful to
embed $\pi_*\oo{_T}$  into the larger sheaf, $\pi_*(j_*\oo{_{T^{reg}}})$.
This is clearly a sheaf of
$\oo{_{T/W}}$-algebras  with a compatible $W$-action along the fibers.

\ab We introduce a sheaf-theoretic analogue of the algebra
 $\C(T)[W]$. It is
a sheaf ${\cal O}[W]$ on $T/W$ (not on $T^{reg}/W$)
defined as the skew tensor product
of $\pi_*(j_*\oo{_{T^{reg}}})$ and $\C[W]$ over $\C$.
Thus, ${\cal O}[W]$ is a sheaf of ${\cal O}_{T/W}$-algebras,
quasi-coherent (but not coherent) as a sheaf of ${\cal O}_{T/W}$-modules.
Clearly, there is
a natural action-morphism of sheaves on $T/W$:
$${\cal O}[W] \otimes \pi_*(j_*\oo{_{T^{reg}}}) \to \pi_*(j_*\oo{_{T^{reg}}})
\eqno (4.1.1)$$

{\bf Definition 4.2.} Let
$\hhh$, resp. $\hk$, be the subsheaf of ${\cal O}[W]$
whose local sections satisfy conditions (1.3.1-2), resp.
(1.3.1-3).
\vskip .1cm

\proclaim Proposition 4.3. Both $\hhh$ and $\hk$ are subsheaves of algebras
in ${\cal O}[W]$. The action (2.1) of $\hhh$ on $\pi_*(j_*\oo{_{T^{reg}}})$
preserves
$\pi_*{\cal O}_{_T}$.

Proof is entirely similar to the proof of Theorem 1.4 and is omitted.

\vskip .2cm
\ab Let ${\J}\subset {\cal O}_{_T}$ be the subsheaf of ideals corresponding
to the divisor
$\sum_{\alpha\in R^{re}_+}$
$T_{\alpha, q^{-2}}$, and $\pi_* {\J} \subset
\pi_*{\cal O}_{_T}$ its direct image.
Our proof of theorem 2.4 yields, in fact,
the following local result.

{\bf Theorem 4.4.} (a) {\it The $\hhh$-action of $\hhh$ on $\pi_*{\cal
O}_T$,
cf. prop. 4.3.,
gives  an isomorphism
$\hhh \simeq End (\pi_*{\cal O}_T)$
of sheaves of algebras on $T/W$. Moreover, we have}
\[\hk= \{ u\in \hhh\;|\; u(\pi_* {\J})\subset\pi_* {\J}\} \qquad\square.\]

\vskip .3cm

\ab Now let
$\cc$ be an elliptic curve. The variety $T/W$ in this case
is known to be a weighted projective space. By the above, the
sheaf $\hk$ is a locally free rank $(\#W)^2$ sheaf of algebras on $T/W$.
We will consider
global sections of $\hk$ over an appropriate open subset.

\ab Let $\xi\in \cc$, $\xi\neq 1$
 be a point of order 2. Observe that the set
$\bigcup_{\alpha\in R^{}_+} T_{\alpha, \xi}$ is
$W$-stable. Its complement is a Zariski-open subset in $T/W$:
\[U=(T\,\setminus\,\bigcup_{\alpha\in R_+} T_{\alpha, \xi})/W\,\]

\ab Choose a complex uniformization $\cc = {\C}/\Lambda$, where
$\Lambda$ is a lattice in $\C$.
Let $sn(z)$ be the Jacobi sine function on ${\C}/\Lambda$ associated with
the choice of the second order point $\xi$. It is uniquely
characterized by the following conditions: it has simple zeroes
at $0$ and $\xi$, simple poles at the two other points of order 2
and its derivative at 0 is equal to 1.
For each simple root $\alpha_i$ we have
$$\sigma_i := {sn(q^{-2})\over  sn(t^{\alpha_i})} [1] +
\biggl (1-{sn(q^{-2})\over  sn(t^{\alpha_i})}\biggl) [s_i]\,\in
\Gamma(U, \hk)
\eqno(4.5)$$
One calculates that, for any $i$, we have
$\sigma_i^2 = [1]$, but these elements do not satisfy the braid relations
[BE] [Gu].

\ab Assume, for instance,  we are in the $SL_2$-case.
Then, $T = \cc$ and
$W = \{1,s\}$. Therefore, we find $T/W=\cc/\{\pm 1\}$ and
$U=(\cc\,\setminus\,\{\xi\})/\{\pm 1\}$.
Let $\wp^{(m)}$ denote the $m$-th derivative of the Weierstrass
$\wp$-function, and write $\sigma$ for the (only) element
in $\Gamma(U, \hk)$ given by formula (4.5).
We have the following result.

{\bf Proposition 4.6.} {\it In the $SL_2$-case
the algebra $\Gamma(U, \hk)$ has $\C$-basis formed by the following elememts:
\[[1]\;,\;
\sigma\;,\; \wp^{(m)}(t-\xi) [1]\;,\;
\bigl(\wp^{(m)}(t-\xi) - \wp^{(m)}(q^{-2}-\xi)\bigl) [s],\]
where $t\in \C/(\{\pm 1\}\ltimes\Lambda)$, and $ m\geq 0$.}

\noindent {\sc Proof:} An element of $\Gamma(U, \hk)$ has the form
$f_1 (t) [1] + f_s(t) [s]$ where $f_1$ and $f_s$ are meromorphic
functions on $\cc$. By definition of the sheaf $\hk$, these functions
are regular for $z\neq 0, \xi$ (where $z$ denotes the coordinate on
$\C$),
  and have at most a simple pole
at $z=0$ with $Res_0 f_1 + Res_0 f_s= 0$.
Set $r = Res_0 f_1$. Then each of the functions
 $f_1 (z) - r/sn(z)$ and  $f_s(z) + r/sn(z)$ is regular off the point
$\xi$ and is therefore a linear combination of $1$ and
$\wp^{(m)}(z-\xi)$. Hence,
we have
$$f_1(t) = c_0 + {c_1\over sn(t)} + \sum_{i\geq 2} c_i \wp^{i-2}
(t-\xi), \; f_s(t) = d_0 - {c_1\over sn(t)} +
\sum_{i\geq 2} d_i \wp^{i-2}(t-\xi)$$
Subtracting linear combinations of 1 and $\sigma$ we can kill
the first two terms in $f_1$ and the second term in $f_s$.
After that the residue condition becomes void, for $\wp^{(i)}$ have
no residues.
Thus, the coefficients
$c_i, i\geq 2$ may be arbitrary. Further, the only remaining condition
on $f_s$ is that $f_s(t)=0$
for $t=q^{-2}$. Clearly, functions
$\wp^{(m)}(t-\xi) - \wp^{(m)}(q^{-2}-\xi)$ form a basis in the space
of all functions satisfying this condition and
 regular off the point $\{t=\xi\}$.
Proposition follows. $\quad\square$

\vskip 7mm
\pagebreak[3]

{\bf 5. Affine Hecke algebra.}

Let $\Phi = (A, \xx)$ be a root
data with Cartan matrix $A$ of finite type. Then,
there is a unique (up to isomorphism) reductive algebraic group $G = G(\Phi)$
over ${\C}$ with Lie algebra ${\g}(A)$ and Cartan subalgebra
 $\h=\C\otimes_{_\Z} \xx^\vee$.

\ab Write $Q\subset {\h}^*$ and  $Q^\vee\subset \h$ for the root and
coroot lattices, respectively, and let $P=Hom_\Z(Q^\vee, \Z)$ and
$P^\vee=Hom_\Z(Q, \Z)$ be the corresponding weight and coweight
lattices. We have
\[Q\subset \xx\subset P \quad\mbox{and}\quad
Q^\vee\subset \xx^\vee\subset P^\vee\]
The group $\Pi_\Phi = \xx^\vee/Q^\vee$ is known to be isomorphic
to the center of $G$.

\ab Let $\theta\in   R_+$ be the maximal root. The affine Weyl group
$ W ^a = W^a(\Phi)$ of $\Phi$ is, by definition, the group
 of affine transformations
of {$\h$} generated by $W$ and $s_0$, additional reflection with respect
to the affine hyperplane $\{ h\in {\h}: \langle\theta, h\rangle+ 1=0\}$.
Thus $W^a$ is a Coxeter group with generators
$s_i, i=0, ..., n$.
It is known that $W^a$ is a semidirect product of $W$ and the co-root
lattice $Q^\vee$.

\ab Let $\widetilde W = \widetilde W(\Phi)$ be the semidirect product of
$W$ and $\xx^\vee\supset Q^\vee$. It is called the extended Weyl group
of $\Phi$ and is not, in general, a Coxeter group. It is clear that
$W^a\subset \widetilde W$ is a normal subgroup and we have
$\widetilde W/W^a = \Pi_\Phi$.
Denote $\widetilde \xx = \xx\oplus {\Z}$, and $\widetilde  R  =
 R\times {\Z} \subset
\widetilde \xx$, the affine root system associated to $R$. Its positive
part is
$$\widetilde R_+ = \bigl\{ (\alpha, k)\in \xx\oplus {\Z}: \alpha\in R,
k>0 \quad {\rm or}\quad \alpha\in  R_+, k\geq 0\bigl\}.$$
The group $\widetilde W$ acts on $\widetilde \xx$ by group automorphisms
preserving $\widetilde R$. Let $\widetilde S$ be the system of simple roots
of $\widetilde  R_+$. It consists of $(\alpha_i, 0)$, $i=1, ..., n$,
which will be still denoted by $\alpha_i$, and $\alpha_0 = (\theta, 1)$.

\ab Write $l(w) := \#{\D}(w)$ for the length of
$w\in W^a\subset \widetilde W$, cf.
sect. 1.5.1. We extend the length function to $\widetilde W$
by the same formula $l(w) := \#{\D}(w)$.
With this definition, the set
$\{w\in \widetilde W: l(w)=0\}$ is  a subgroup which is
identified with  $\Pi_\Phi$. The group $\Pi_\Phi$ acts on $\widetilde S$
by permutations. Thus $\widetilde W = \Pi_\Phi \ltimes W^a$,
the semidirect product with the commutation relations
$\pi s_\alpha \pi^{-1} = s_{\pi(\alpha)}$, $\alpha\in \widetilde S$.

\ab Let $q\in {\C}^*$ be a non-zero
complex number. The affine Hecke algebra $\hd_q = \hd_q(\Phi)$
can be defined in several different ways.

\noindent {\bf (5.2) First definition:}
 $\hd_q$
is the ${\C}$-algebra with basis $T_w, \, w\in \widetilde W$ and
multiplication given by the rules:
$$(T_w+q^{-1})(T_w - q) = 0, \quad w\in \{s_0, ..., s_n\},\eqno (5.2.1)$$
$$T_wT_y = T_{wy}, \quad {\rm if}\quad l(wy) = l(w) + l(y).
\eqno (5.2.2)$$
Denote by $\hd_q^a\subset \hd_q$ the subspace spanned by the $T_w$,
$w\in W^a$ only. This is clearly a subalgebra and $\hd_q \simeq
 \hd_q^a[\Pi_\Phi]$, is the twisted group algebra for the $\Pi_\Phi$-action
on $\hd_q^a$. In other words, $\hd_q$ is generated by the sets
$\{T_w, w\in W^a\}$ and  $\{T_\pi, \pi\in \Pi_\Phi\}$  with  the relations
(5.2.1-2) for the  $T_w$'s, and the relations
$$T_\pi T_{\pi'} = T_{\pi+\pi'}; \,\,\,
T_\pi T_{s_\alpha}T_\pi = T_{s_{\pi(\alpha)}},
\alpha \in \widetilde
R_+.\eqno (5.2.3)$$

\noindent {\bf (5.3) Second definition:} We have
$\widetilde W\simeq W \ltimes \xx^\vee$.
Accordingly, one has a presentation of $\hd_q$ by generators $\{T_w,
w\in W\}$ and $\{Y_\lambda, \lambda \in \xx^\vee\}$ subject to the relations:

(5.3.1) The $T_w, w\in W$, satisfy (4.2.1-2).

(5.3.2) $Y_\lambda Y_\mu = Y_{\lambda+\mu}.$

(5.3.3) $T_{s_i} Y_\lambda = Y_\lambda T_{s_i}$,
 if $(\alpha_i^\vee, \lambda) = 0$.

(5.3.4) $T_{s_i} Y_{s_i(\lambda)}T_{s_i} = q Y_\lambda$, if
$(\alpha_i^\vee, \lambda) = 1$.

It is known that the elements $T_w Y_\lambda$, $w\in W$, $\lambda\in \xx$,
form a ${\C}$-basis of $\hd_q$.
\vskip 2mm
\ab Recall the Demazure-Lusztig elements $\sigma_i\in \Hk(\Phi)$
 given, for $i=1, ..., n$, by formula (2.4).
Theorem 1.7 is a consequence of the following
more precise result.
\vskip .2cm

{\bf Theorem 5.4.} {\it The assignment $T_{s_i}\mapsto \sigma_i$,
$Y_\lambda\mapsto t^\lambda$ extends to an algebra isomorphism
$f: \hd_q\rightarrow \Hk(\Phi^\vee)$, where $\Phi^\vee$
is the root data dual to  $\Phi$.}

\noindent {\sc Proof:} The existence of $f$ is proved by verifying the
relations, which is well known. The fact that $f$ is  an isomorphism
follows since $\{T_w Y_\lambda\}$ is a ${\C}$-basis in $\hd$, while
$\{\sigma_w t^\lambda\}$ is a ${\C}$-basis in $\Hk$, by theorem 2.5.

\vskip .3cm
\pagebreak[3]

{\bf 6. Cherednik algebra.}

{\bf (6.1)} {\sc Affine root systems reviewed.}
Let $A$ be a generalized Cartan matrix of untwisted affine
type. We write it
as $A = \|a_{ij}\|_{i,j=0}^n$, so that:

${\overline A}= \|a_{ij}\|_{i,j=1}^n$ is a (finite type) Cartan matrix
from which $A$ is obtained by affinization.

${\overline \g}= {\g}({\overline A})$ is the corresponding finite-dimensional
Lie algebra.

${\overline G}$  is the simply connected reductive group over ${\C}$ with
Lie algebra ${\overline \g}$.

${\overline R}, {\overline W}, {\overline T}$
etc., are the root system, the Weyl group,
the maximal torus, etc. for ${\overline G}$. Thus, $R = \{(\alpha, k),
\alpha\in {\overline R}, \, k\in {\Z}\}$. The simple root $\alpha_0$
is $(-\theta,1)$ where $\theta$ is the maximal root of ${\overline R}$.
Notice that the lattice of characters
of ${\overline T}$ is ${\overline Q}$, the dual of the coroot lattice of
${\overline \g}$.

$\Pi = {\overline P}^{\vee}/{\overline Q}^{\vee}$ is the center of
the simply connected group associated to ${\overline\g}$.

\noindent We have the identification of vector spaces :

${\g} = \g(A) = {\C}\cdot d \oplus {\C}\cdot c
\oplus {\overline \g}[z, z^{-1}]$, where $c$ is the central element, $d$ is
the infinitesimal rotation, i.e., $[d, g(z)] = z\frac{dg}{dz}$ for
$g(z)\in {\overline \g}[z, z^{-1}]$.

\vskip .1cm

${\g}_{der} = \g_{der}(A)= {\C}\cdot c
\oplus {\overline \g}[z, z^{-1}]$ is the derived algebra of $\g$,
considered  for the purpose of comparing
our construction with that of Cherednik.

\vskip .1cm

$G_{der} = \widehat{{\overline G}((z))}$ is the Kac-Moody group associated with
${\g}_{der}$. Its center will be denoted ${\C}^*_\zeta$.
It is isomorphic to ${\C}^*$ and we think of $\zeta$ as a coordinate
therein.

\vskip.1cm

${\h}_{der} = {\C}\cdot c \oplus {\overline\h}$  is the Cartan subalgebra
of ${\g}_{der}$, and
$T_{der}  = {\C}^*_\zeta \times {\overline T}$ the maximal torus of
$G_{der}$.

\ab The $W$-action on ${\C}[T_{der}]$, the
ring of regular functions on  $T_{der}$, is given by the formulas:
$$s_i (t^\lambda \zeta^l) = t^{s_i(\lambda)}\zeta^l, \quad i=1, ..., n,
\, \lambda\in {\overline P}, l\in {\Z},\eqno (6.1.1)$$
$$s_0 (t^\lambda \zeta^l) = t^{\lambda-<\lambda,\theta^\vee>\theta}\zeta^{l+1}.
\eqno (6.1.2)$$

{\bf (6.2)} The  double affine Hecke algebra
$\hdd_q$ associated to $A$ was defined by Cherednik [Ch2] to be the
$\C$-algebra
(depending on a parameter $q\in {\C}^*$, $q\neq \pm 1$) on generators
$$\delta\;,\quad T_i\,,\, i=0, ..., n\;,\;
\quad    Y_\lambda\,,\,
\lambda\in\overline P\;,\quad T_\pi\,, \,\pi\in \Pi,$$
subject to the following relations:

(6.2.0) The element $\delta$ is central.

(6.2.1) The $T_i, T_\pi$ generate the affine Hecke algebra $\hd_q$
corresponding to the finite root system $(\overline A,\overline P)$.

(6.2.2) The $T_i, i=1, ..., n$ and the $Y_\lambda$ generate the
affine Hecke algebra corresponding to the dual
finite root system $(\overline A^t,\overline P^\vee)$.

(6.2.3)
$T_i Y_\lambda T_i = Y_\lambda Y_{\alpha_i}^{-1}$,
if $<\lambda, \alpha_i^\vee>=1$, $i=1, ..., n$.

(6.2.4) $T_0^{-1} Y_\lambda T_0^{-1} = Y_\lambda Y_\theta^{-1}\delta$,
if $<\lambda, \theta^\vee>=1$.

(6.2.5) $T_i Y_\lambda = Y_\lambda T_i$ for any $i=0,..., n$,
$\lambda\in \overline P$
such that $<\lambda, \alpha_i^\vee>=0$. Here we
set $\alpha_0^\vee = -\theta^\vee$.

(6.2.6) $T_\pi Y_\lambda T_\pi^{-1} = Y_{\pi(\lambda)}$.

\vskip .2cm

\ab We are going to compare
the algebra $\hdd_q$ with ${\bold H}_q(\Phi_{der})$, where $\Phi_{der}=(A,\xx)$
is the root data such that $\xx\subset\h_{der}^*$ is the
lattice dual to $\xx^\vee = \bigoplus_{i=0}^n\alpha_i^\vee\subset\h_{der}$.
The subalgebra in $\hdd_q$
generated by $\delta, Y_\lambda$, $\lambda\in\overline P$ is naturally
identified with ${\C}[T_{der}]$, where $\delta$ corresponds to the
coordinate function, $\zeta$, on the center of the derived Kac-Moody group.
Recall also the elements $\sigma_i, i=0, ..., n$ of ${\bold H}_q(\Phi_{der})$
from (2.4).
\vskip .2cm

{\bf Theorem 6.3.1.} {\it The assignment
$\;Y_\lambda\delta^l\mapsto t^\lambda\zeta^l$ , $ T_i \mapsto \sigma_i\,,\,
i=0, ..., n\;$
gives an algebra isomorphism} ${\hdd}_q\rightarrow {\bold H}_q(\Phi_{der})$.

\noindent {\sc Proof.} We have first to verify that the relations
of $\hdd_q$ hold in ${\bold H}_q(\Phi_{der})$, which is known, see, e.g., [Ki].
This means that we have a homomorphism, denote it $f$.
To show that it is an isomorphism
we use the fact, proved by Cherednik, that the elements
$Y_\lambda\delta^l
T_w$, $w\in W$, form a ${\C}$-basis in ${\hdd}_q$. A similar statement
holds for ${\bold H}_q(\Phi_{der})$, by Theorem 2.5. This implies that
$f$ gives isomorphisms on the factors of the natural filtrations in
${\bold H}_q(\Phi_{der})$ and $\hdd_q$, labeled by $W$, and thus is an
isomorphism.

\vskip 1cm
\pagebreak[3]
\centerline {\bf References}

{\small

[Bou] $\hspace{3mm}$ \parbox[t]{115mm}{
N. Bourbaki, {\it
Groupes et Alg\`ebres de Lie.} Masson, Paris, 1981.}

[BE] $\hspace{3mm}$ \parbox[t]{115mm}{
 P. Bressler, S. Evens, The Schubert calculus, braid relations and generalized
cohomology, {\it Trans. AMS} {\bf 317} (1990), 799-811.}

[Ch1] $\hspace{2mm}$ \parbox[t]{115mm}{
I.V. Cherednik,
Double Affine Hecke
algebras, KZ-equations and Macdonald's operators.
{\it International Mathem. Research Notices, Duke Math. J.},
{\bf 9} (1992), 171-180.}

[Ch2] $\hspace{2mm}$ \parbox[t]{115mm}{
I.V. Cherednik, Double affine Hecke algebras and Macdonald conjectures,
{\it Ann. Math.} 1995.}

[Dr] $\hspace{3mm}$ \parbox[t]{115mm}{
V. Drinfeld, Degenerate affine Hecke algebras and Yangians,
{\it Funct. Anal. and Appl.} {\bf 20} (1986), 69-70.}

[GG] $\hspace{2mm}$ \parbox[t]{115mm}{
H. Garland, I. Grojnowski, "Affine" Hecke algebras associated to
Kac-Moody
groups. Yale preprint 1995 (q-alg/9508019).}

[Gu] $\hspace{2mm}$ \parbox[t]{115mm}{
E. Gutkin, Representations of Hecke algebras,
{\it Trans. AMS}, {\bf 309} (1988),
269-277. }

[Ha] $\hspace{3mm}$ \parbox[t]{115mm}{
R. Hartshorne, Algebraic Geometry, Springer-Verlag 1983.}

[Ka] $\hspace{3mm}$ \parbox[t]{115mm}{V.G. Kac, Infinite-dimensional Lie
algebras,
 Cambridge University Press, 1991.}

[KL] $\hspace{3mm}$\parbox[t]{115mm}{
D. Kazhdan, G. Lusztig, Proof of the Deligne - Langlands conjecture for Hecke
algebras, {\it Invent. Math.} {\bf 87} (1987), 153-215.}

[Ki] $\hspace{3mm}$ \parbox[t]{115mm}{
A. Kirillov, Jr. Lectures on the affine Hecke algebras and Macdonald
 conjectures, Harvard preprint, 1995. }

[KK1] $\hspace{1mm}$ \parbox[t]{115mm}{
B. Kostant, S. Kumar. The Nil-Hecke ring and the cohomology
of $G/P$ for a Kac-Moody group $G$. {\it Advances Math.} {\bf 62}
(1986),
187-237.}

[KK2] $\hspace{1mm}$ \parbox[t]{115mm}{
B. Kostant, S. Kumar. $T$-equivariant $K$-theory of generalized
flag varieties. {\it J. Diff. Geom.} {\bf 32} (1990), 549-603.}

[L1] $\hspace{3mm}$ \parbox[t]{115mm}{
G. Lusztig. Equivariant $K$-theory and representations of affine Hecke
algebras, {\it Proc. of the A.M.S.} {\bf 94} (1985), 337-342.}

[L2] $\hspace{3mm}$ \parbox[t]{115mm}{
G. Lusztig. Cuspidal local systems and graded Hecke algebras.
{\it Publ. Mathem. I.H.E.S.} {\bf 67} (1988), 145-212.}

[Sp] $\hspace{3mm}$ \parbox[t]{115mm}{
T. Springer, Reductive groups, {\it Proc. Symp. Pure Math.} {\bf 33},
part 1, p. 3-27, Amer. Math. Soc., Providence RI 1979.
\vskip 1cm
\nopagebreak{
V.G.: Department of Mathematics, University of Chicago, Chicago IL
60637, USA;
ginzburg@@math.uchicago.edu

\vskip .3cm

M.K.: Department of Mathematics, Northwestern University,
Evanston IL 60208, USA;
kapranov@@chow.math.nwu.edu

\vskip .3cm

E.V.: Universit\'e de Cergy-Pontoise, 2 av. A. Chauvin, 95302 Cergy-Pontoise
Cedex, France;
vasserot@@paris.u-cergy.fr}

}
\end{document}